\documentstyle{amsppt}
\magnification=\magstep1 \NoBlackBoxes \TagsOnRight
%\pageno=47
\baselineskip 20pt \pagewidth {5 in} \pageheight {6.75 in}
\NoRunningHeads \topmatter
\title Cup-length estimate for Lagrangian intersections
\endtitle
\author Chun-gen Liu \endauthor
\affil Department of Mathematics, Nankai University\\
Tianjin 300071, People's Republic of China \\
liucg\@nankai.edu.cn\endaffil
\thanks
* Partially supported by the National Natural Science Foundation of
China(10071040), Ph.D Fund of ME of China, PMC Key Lab of ME of
China.
\endthanks
\abstract In this paper we consider the Arnold conjecture on the
Lagrangian intersections of some closed Lagrangian submanifold of
a closed symplectic manifold with its image of a Hamiltonian
diffeomorphism. We prove that if the Hofer's symplectic energy of
the Hamiltonian diffeomorphism is less than a topology number
defined by the Lagrangian submanifold, then the Arnold conjecture
is true in the degenerated (non-transversal) case.
\endabstract
\endtopmatter
\def\<{\langle}
\def\>{\rangle}
\def\r{\Bbb R}

\def\n{\Bbb N}
\def\z{\Bbb Z}

\def\u{\nu_{\tau}}

\def\u{\tilde {u}}
\def\M{\Cal M(J,H,x^+,x^-)}
\def\m{\tilde{\Cal M}(J,H,x^+,x^-)}
\document
\head \S1 Introduction and main results\endhead

Let $(M,\omega)$ be a closed symplectic manifold, $L\subset M$ be
its closed Lagrangian submanifold. A Hamiltonian $H:[0,1]\times
M\to \r$ is a $C^{\infty}$ function. This function defines a
$t$-dependent Hamiltonian vector field $X_{H_t}$ on $M$ by
$\omega(\cdot,X_{H_t})=dH_t$. The time one map $\varphi=\varphi^1$
of the flow generated by the Hamiltonian vector field $X_{H_t}$ is
a symplectic automorphism of $M$. Arnold conjecture that, for some
symplectic manifold $(M,\omega)$ and its Lagrangian submanifold
$L$, the intersection $L\cap \varphi(L)$ contains at least as many
points as a topology number of $L$. If $L$ transversely meet
$\varphi(L)$, then the topology number can be the
$\text{rank\;of\;}H^*(L;\Bbb F)$ for some ring or field $\Bbb F$.
In general, this topology number can be the cup-length of $L$
which is defined by
$$ \aligned \text{cl}(L,\Bbb F)=\max\{k+1|&\;\exists\alpha_i\in H^{d_i}(L,\Bbb F),\; d_i\ge 1,\;
i=1,\cdots,k\;\\&\text{such that}\;
\alpha_1\cup\cdots\cup\alpha_k\ne 0 \}.\endaligned$$
 In this paper, we fixed $\Bbb F=\z_2$ and denote the cup-length of $L$ by $\text{cl}(L)$.

It is well known that the above Arnold conjecture is not true in
general. For example the ``small Lagrangian torus" in a symplectic
manifold can be push away by some Hamiltonian diffeomorphism. In
this case the intersection $L\cap \varphi(L)=\varnothing$, but the
topology number of $L$ is not zero. So we need further conditions
to guarantee this version of the Arnold conjecture. The first
condition was given by Floer in [F1,F2] (see also [H]). It was
proved that if $\pi_2(M,L)=0$ or $\omega(\pi_2(M,L))=0$, then the
Arnold conjecture on the Lagrangian intersection is true. Chekanov
[Ch1-Ch2] found that there is some relation between the Hofer's
bi-invariant metric of the Hamiltonian diffeomorphim and this
version of Arnold conjecture.

For a Hamiltonian $H:[0,1]\times M\to \r$, we can define a
semi-norm of $H$ as
$$\|H\|=\int^1_0(\max_x\,H(t,x)-\min_x\,H(t,x))\,dt. $$
This semi-norm is weaker than $C^{0}$-norm of $H$ and plays an
eminent role for Hofer's bi-invariant metric on the group of
compactly supported Hamiltonian diffeomorphism. The metric is
defined by
$$d(\varphi,id_M)=\inf\{\|H\|\;|\,\varphi\,\text{is\, generated\, by\, }H\}.$$

 We say that $L$ is a rational Lagrangian submanifold of $M$ if
there is a number $\sigma(L)>0$ such that
$\omega(\pi_2(M,L))=\sigma(L)\cdot \z$.

 Chekanov in [Ch1] (see [Ch2] for a somewhat general statement)
 proved that if $d(\varphi,id_M)<\sigma(L)$, then $\sharp (L\cap \varphi(L))\ge \dim
 H^*(L;\z_2)$ provided $L$ is a rational Lagrangian submanifold of
 $M$ with the number $\sigma(L)>0$ defined as above and the
 intersection is transverse.

 The main result of this paper is the following theorem.
 \proclaim{Theorem} If $L\subset M$ is a rational Lagrangian submanifold of
 $M$ with the number $\sigma(L)$ defined above and
 $d(\varphi,id_M)<\sigma(L)$, then there holds
 $$\sharp (L\cap \varphi(L))\ge  \text{\rm cl}(L).$$
  \endproclaim

{\bf Acknowlegements} Part of this  work  was completed when the
author was visiting Max Planck Institute for Mathematical Science
in Leipzig, Germany. The author would like to thank J. Jost, M.
Schwarz, G. Wang, Jost Li and A. Vandenberg for their hospitality,
many helps and helpful conversations.

 \head \S2 $J$-holomorphic curves with boundary conditions
\endhead

Let $L$ be a closed embedded Lagrangian submanifold of a compact
symplectic manifold $(M,\omega)$. $H: [0,1]]\times M\to \r$ is a
smooth function, and $\varphi^t$ is the Hamiltonian flow generated
by the Hamiltonian function $H$. Setting $L_1=\varphi^1(L)$, and
considering the space
 $$ \Omega_1(L)=\{\gamma\in C^{\infty}([0,1],M)\,|\,\gamma(0)\in L,\;\gamma(1)\in L_1\},$$
 restricting to this space we define a $1$-form $\alpha$ by
 $$\<\alpha(\gamma),\xi\>=\int^1_0\omega(\dot \gamma(t),\xi(t))\,dt. $$
 This $1$-form is closed. Let $\Omega^0_1(L)$ be the component of
 $\Omega_1(L)$ which contains the constant path. A primitive $F$
 of $\alpha|_{\Omega^0_1(L)}$ is a $\r/\sigma\z$-valued functional
 on $\Omega^0_1(L)$, the standard action functional of Floer's
 theory. It is defined up to additive constants. For a compatible almost complex
structure $J$, define a metric on $\Omega_1(L)$ as follows:
$$\<\xi_1,\xi_2\>=\int_0^1\omega(\xi_1(t),J\xi_2(t))\,dt. $$
The gradient of $F$ with respect to this metric is given by
$$\nabla F(\gamma)(t)=J(\gamma(t))\dot \gamma(t).$$
For a pair $(x^+,x^-)$ of critical points of $F$ which
correspondent to a pair of intersections of $L\cap L_1$, we
consider the following moduli space which is analogue to the
connect orbit space of the negative gradient flow of a Morse
functional defined on a finite dimensional space
$$\Cal M(J,H,x^+,x^-)=\left\{u:\r\times [0,1]\to M \Big |   \,
\aligned & \partial_s u+J\partial_t u=0,\;u\;\text{is\,not\,constant},\\
& u(s,0)\in L, \;\; u(s,1)\in \varphi^{1}(L)\\
& \lim_{s\to \pm \infty}u(s,t)=x^{\pm}\in L\cap
\varphi^1(L)\endaligned \right \}$$
 If $u\in \Cal M(J,H,x^+,x^-)$, we define a map $\tilde{u}:\r\times [0,1]\to M$
such that $u(s,t)=\varphi^t(\u(s,t))$, then we get
$$ \partial_s \tilde{u}+\tilde {J}_t\left(\partial_t \tilde u
+X_{\tilde{H}}(\tilde {u}(s,t))\right)=0.\tag 2.1$$
Here $\tilde{J}_t=(d\varphi^t)^{-1}Jd\varphi^t$,
$\tilde{H}(t,x)=H(t,\varphi^t(x))$ and
$X_{\tilde{H}}(x)=(d\varphi^t)^{-1}X_H(\varphi^t(x))$ by
definition. If $J$ is compatible with the symplectic structure
$\omega$, so is for the $t$-dependent almost complex structure
$\tilde {J}_t$. $\tilde u$ satisfies the following conditions
(2.2) and (2.3).
$$\cases & \u(s,0)\in L,\;\; \forall s\in (-\infty,+\infty)\\
         & \u(s,1)\in L, \;\;\forall s\in (-\infty,+\infty). \endcases \tag 2.2  $$
$$\lim_{s\to \pm \infty}\u(s,t)=x^{\pm}(t), \tag 2.3 $$
where $x^{\pm}(t)=(\varphi^t)^{-1}(x^{\pm})$ is a Hamiltonian
flow line of the Hamiltonian function $-\tilde H$ and
$x^{\pm}(0)=x^{\pm}\in L\cap \varphi^1(L)$. Conversely, if $\u$
is a solution of (2.1) satisfies (2.2) and (2.3), then
$u(s,t)=\varphi^t(\u(s,t))$ belongs to $\M$. In fact, it is easy
to see $u$ solves the equation
$$ \partial_s u+J \partial_t u
 =0. \tag 2.4 $$
 By definition of $u$, we have
 $$u(s,0)=\u(s,0)\in L,\;\;u(s,1)=\varphi^1(\u(s,1))\in \varphi^1(L), \tag 2.5  $$
  and
  $$\lim_{s\to\pm\infty}u(s,t)=\varphi^t(x^{\pm}(t))=x^{\pm}(0)\in L\cap \varphi^1(L). \tag 2.6 $$
  Thus we can consider the following moduli space
  $$\m=\left\{\u \Big |   \, \aligned & \partial_s \u  +\tilde{J}_t\left( \partial_t \u
  +X_{\tilde{H}}(\tilde {u}(s,t))\right)=0\\
& \u(s,0)\in L, \;\; \u(s,1)\in  L\\
& \lim_{s\to \pm \infty}\u(s,t)=x^{\pm}(t)\, \text{is Hamiltonian
flow line of}\, -\tilde{H},\\ & \;\;x^{\pm}(0)=x^{\pm}\in L\cap
\varphi^1(L)\endaligned \right\}.$$
 This moduli space $\m$ is 1-1 corespondent with  $\M$.

 We recall that the Hamiltonian flow line of the Hamiltonian
 function $-\tilde H$ with Lagrangian boundary condition is a
 solution of the following equation
 $$\cases &\dot x(t)=-X_{\tilde H}(x(t))\\& x(0)\in L,\;\;x(1)\in L.\endcases \tag 2.7$$
 We can write $x(t)=(\varphi^t)^{-1}(x_0)$, then $x(0)=x_0\in L$ and
 $x(1)=(\varphi^1)^{-1}(x_0)\in L$, it implies $x(0)=x_0\in
 L\cap\varphi^1(L)$. The space of the solutions of (2.7) is one to
 one correspondent with the set $L\cap \varphi^1(L)$.

 In order to find  solutions of equation (2.7), we define the
 following spaces
 $$\aligned &\tilde \Omega(L)=\{ x\in C^{\infty}([0,1],M)\;|\,\;x(0)\in L,\;\;x(1)\in L
 \},\\
 &\tilde \Omega_0(L)=\{ x\in \Omega(L)\;|\;\,[x]=0\in \pi_1(M,L) \}, \endaligned$$
 and the universal cover space of  $\tilde \Omega_0(L)$
 $$\Omega_0(L)=\{ u_x:D\to M\;|\;u_x|_{S^+}=x , \;u_x|_{S^-}=\tilde x\},$$
where $D$ is the unit disc in $\Bbb C$ with $\partial D=S^+\cup
S^-$, and $S^+$ (resp. $S^-$) is the upper (resp. lower) half unit
circle which is a part of $\partial D$, the boundary of $D$.
$\tilde x:[0,1]\to L$ is a path in $L$ which is isotopic to $x$
relative to the end points. On the space $\Omega_0(L)$ we define a
functional
$$\Cal A_H(x,u_x)=\int_Du_x^*\omega+\int_0^1\tilde H(t,x(t))\,dt. $$
It is easy to see that
$$d\Cal A_H(x)(\xi)=\int_0^1\omega(\dot x+X_{\tilde H}(x),\xi) $$
This means that $d\Cal A_H(x)=0$ implies $\dot x+X_{\tilde
H}(x)=0$.

The functional induces a functional $\tilde A:
\tilde\Omega_0(L)\to \r/\sigma\z$ if $L$ is rational with
$\omega({\pi_2(M,L)})=\sigma(L)\z$ for some $\sigma=\sigma(L)>0$.

\head\S3 Morse homology and its cup product\endhead

We first recall the Morse homology theory briefly
 (see [MS2] for details), Let $(f,g)$ be a Morse-Smale pair on
 $L$, that is, let $f$ be a fixed Morse function and $g$ be a
 generic Riemannian metric on $L$ such that the stable and
 unstable manifolds $W^s(y)$, $W^u(x)$ for critical points  $x,y\in \text{Crit}
 f$ for the negative gradient flow of $(f,g)$ intersect
 transversely.  We define the connect orbit space of $x,y\in
 \text{Crit}f$ by
 $$M_{x,y} (f,g)=\{ \gamma\in C^{\infty}(\r,L)\,|\;\dot \gamma+
 \nabla_{g}f(\gamma)=0, \;\gamma(-\infty)=x,\;
 \gamma(+\infty)=y \}.$$
 We have $\dim M_{x,y}(f,g)=\mu(x)-\mu(y)$, $\mu(x)$ is the Morse index of $x\in  \text{Crit}f$,
  and $M_{x,y}(f,g)$
 admits a free $\r$-action by translation: $s\cdot \gamma(\cdot)=\gamma(s+\cdot)$.
 We denote the quotient space by
 $$\hat M_{x,y}(f,g)=M_{x,y}(f,g)/\r. $$
Let $C^k(f)$ denote the $\z_2$-free Abelian  group generated by $
\text{Crit}_k f=\mu^{-1}(k)$, and define the boundary operator as
$$\delta:C^k(f)\to C^{k+1}(f),\;\;\delta x=\sum_{\mu(y)=\mu(x)+1}n(y,x)y $$
where $n(x,y)$ is defined by
$$n(x,y)=\sharp _{\z_2} \hat M_{x,y}(f,g) $$
the modulo 2 number of $\hat M_{x,y}(f,g)$, it is well defined
when $\mu(x)-\mu(y)=1$. It is well known that $\delta^2=0$, and
$$H^*(C^*(f),\delta)\cong H^*(L;\z_2). \tag 3.1$$
Let $(f,g_i),\;i=1,2,3$ be three generic Morse-Smale pairs on $L$
such that the following moduli spaces are $\mu(z)-\mu(x)-\mu(y)$
dimensional space for $x,y,z\in {\text Crit}f$
$$\aligned \Cal M_{z,x,y}(f,g_1,g_2,g_3)=\{(\gamma_1,\gamma_2,\gamma_3)
&\in W^u(z)\times W^s(x)\times W^s(y)\;|\\
&\gamma_1(0)=\gamma_2(0)=\gamma_3(0)\}\endaligned $$ and the
spaces $\Cal M_{z,x,y}(f,g_1,g_2,g_3)$ are compact in dimension
$0$.

Analogously to $\delta$ we define the following operation on
$C^*(f,\z_2)$. Given $x,y,z\in {\text Crit}f$, we set
 $$n(z;x,y)=\sharp \Cal M_{z,x,y}(f,g_1,g_2,g_3)\;(\text{mod}\; 2)
 \;\text{for}\;\mu(z)=\mu(x)+\mu(y) $$
 and
 $$m_2:C^k(f,\z_2)\otimes C^l(f,\z_2)\to C^{k+l}(f,\z_2) $$
 $$m_2(x\otimes y)=\sum_{z\in {\text Crit}_{k+l}f}n(z;x,y)z. \tag 3.2$$
 $m_2$ is a chain operator and it induced a cup product of the
 cohomologies $H^*(L;\z_2)$.
 These result are standard now (see for example: [MS1] section 3 for $A=0$
 thus $u$ must be a constant map, or [Fu1] for
 $f_1=f_2=f_3$ with different metrics satisfying the transversal
 conditions). Analogously we can define the moduli spaces for $x_0, x_1, \cdots, x_k\in {\text Crit}f$
 $$\aligned \Cal M_{x_0; x_1, \cdots, x_k}=\{(\gamma_0,\gamma_1,\cdots,\gamma_k)&\in W^u(x_0)\times W^s(x_1)
 \times\cdots\times W^s(x_k)\;|\\
&\gamma_0(0)=\gamma_1(0)=\cdots=\gamma_k(0)\}\endaligned  $$ and
$$m_k: C^{l_1}(f,\z_2)\otimes\cdots\otimes C^{l_k}(f,\z_2)\to C^{l_1+\cdots+l_k}(f,\z_2) $$
$$m_k(x_1,\cdots,x_k)=\sum_{x_0}n_k(x_0;x_1,\cdots,x_k)x_0,\; where$$
$$\mu(x_0)=\mu(x_1)+\cdots+\mu(x_k)\;\;\text{and}\;\; n_k(x_0;x_1,\cdots,x_k)=
\sharp_{\z_2} \Cal M_{x_0; x_1, \cdots, x_k}.$$
 $m_k$ induced
$k$-fold cup-product of the cohomologies $H^*(L;\z_2)$.

In this section we always assume that $(M,\omega)$ is a closed
 symplectic manifold. $L\subset M$ is a closed rational Lagrangian submanifold with the
 constant $\sigma(L)>0$ defined as in section 2. i.e., we have
 $\omega({\pi_2(M,L)})=\sigma(L)\z $ for some $\sigma(L)>0. $
 Denote by $\Cal H(M)$ the set of all Hamiltonian function $H:[0,1]\times M\to \r$ .
 Any $H\in \Cal H(M)$ defines a time-dependent  Hamiltonian flow $\varphi^t:M\to M$.
 Time one maps of such flows form a group $\Cal S(M,\omega)$
 called the group of Hamiltonian symplectomorphisms of $M$.
 On the space $\Cal H(M)$, we have a semi-normal defined by
 $$\|H\|=\int^1_0(\max_xH(t,x)-\min_xH(t,x))\,dt. $$
 For $\varphi\in \Cal S(M,\omega)$, the energy of $\varphi$ is
 defined by
 $$E(\varphi)=\inf\{ \|H\|\,|\,\varphi \text{ is  a  time  one  flow  generated  by  }
   H\in \Cal H(M)\} $$
We assume that $E(\varphi)<\sigma(L)$, this condition is essential
for the compactness of the moduli spaces because under this
condition no bubbling-off ($J$-holomorphic sphere and disc)
occurs. So we can naturally define the deformation cup product of
the cohomology groups. Under the above conditions, we have the
moduli space
$$\Cal M^0(J,H)=\{\u\in C^{\infty}(D,M)\;|\;  \partial_s \u
 +\tilde{J}_t\left( \partial_t \u
 +X_{\tilde{H}}(\tilde {u}(s,t))\right)=0,\;[\u]=0\in
\pi_2(M,L)\}.$$
 Given  $x_0, x_1, \cdots, x_k\in {\text Crit}f$ we
define
 $$\aligned \Cal M^0_{x_0; x_1, \cdots, x_k}=\{&(\u,\gamma_0,\gamma_1,
 \cdots,\gamma_k)\in \Cal M^0(J,H)\times W^u(x_0)\times W^s(x_1)
 \times\cdots\times W^s(x_k)\;|\\
&\u(z_i)=\gamma_i(0),\;z_i\in
\partial D,\;i=0,1,\cdots,k\}\endaligned$$
 \proclaim{Theorem 3.1} Given
a   Hamiltonian function $H$ with $\|H\|<\sigma(L)$, and generic
pairs $(f,g_i),\;i=0,1,\cdots,k$, the following operator $m^0(H)$
is well defined,
 $$ m_k^0(H):C^{l_1}(f)\otimes \cdots \otimes C^{l_k}(f)\to C^{l_1+\cdots +l_k}(f),$$
 $$m_k^0(H)(x_1\otimes\cdots\otimes x_k)=\sum_{x_0} (\sharp \Cal M^0_{x_0; x_1, \cdots, x_k}\;\text{mod}\;2)x_0$$
Moreover, $m^0(H)$ is a co-chain map with respect to the boundary
operator $\delta$, and the induced  operation of the cohomology
group is just the  $k$-fold cup product in the sense of (3.1).
\endproclaim
\demo{Proof} The essential ingredient of the proof is to prove the
fact of no bubbling-off. This can be done by looking at the energy
of the element $\u\in\Cal M^0(J,H)$

$$\aligned E(\u)=&\int_{D}|\partial_s \u |_{\tilde J}^2\,dsdt=
\int_D\omega(\partial_s \u, \tilde{J}_t\partial_s \u)\\
 =&-\int_D \omega(\tilde{J}_t\left({\partial_t \u}
 +X_{\tilde{H}}(\tilde {u}(s,t)),\tilde{J}_t\partial_s \u)\right)\,dsdt\\
 =&- \int_D \omega(X_{\tilde{H}}(\tilde {u}(s,t), \partial_s \u)\,dsdt \\
 =& \int_D  d\tilde {H_t}(\u(s,t))(\partial_s \u) dsdt\le \|H\|.\endaligned\tag 3.3$$
 Here we have use the condition $[\u]=0\in \pi_2(M,L)$. Since
 $\|H\|<\sigma(L)$,  notice that we can take $\pi_2(M)$ as a
 sub-group of $\pi_2(M,L)$, any bubbling-off must have energy at
 least $\sigma(L)$, so no bubbling-off occurs. If $H\equiv 0$,
 then we have $m^0(0)=m_k$ as defined in (3.2) which induced the $k$-fold cup
 product. Taking a suitable homotopy $H\sim 0$ such that the
 induced maps in $H^*(L; \z_2)$ satisfying
$m_k^0(H)^*=m_k^0(0)^*=m_k^*$ (see [MS1], Theorem 3.8 for similar
arguments. Here we only consider $A=0$).
 \hfill $\blacksquare$ \enddemo

 \head \S4 The proof of the main result\endhead
 We follow the ideas of [MS1] to prove the main result of this
 paper.
 Firstly, we modify the pair $(J, H)$ and define the ``adapted solution spaces".
  Given the Hamiltonian $H\in C^{\infty}([0,1]\times M,\r)$
 and an $\omega$-compatible almost structure $J$, we get a
 corresponding pair $(\tilde J, \tilde H)$ as in section 2. Here
 $\tilde J$ is explicitly dependent of $t\in [0,1]$. Pick an
 $t$-independent almost complex structure $J_0$ on $TM\to M$, we
 extend $\tilde J$ and $J_0$ to a smooth $1$-parameter family $\bar J=\bar
 J(s)$, $s\in (-\infty,+\infty)$ as
 $$\bar J(s)=\cases J_0,\;\;& s\le 0,\\
 \tilde J,\;\; & s\ge 1.\endcases  \tag 4.1$$
 Let $\beta\in C^{\infty}(\r,[0,1])$ be a monotone cut-off
 function such that
 $$\beta(s)=\cases0,\;&s\le 0,\\
 1,\;&s\ge 1, \endcases \; \; \text{and} \; \beta'(s)\ge 0.$$
 For $R\in [1,\infty)$, we  defined $1$-parameter pairs $(\tilde J_R,\tilde
 H_R)$ on $\r\times [0,1]\times M$ as follows,
 $$(\tilde J_R,\tilde H_R)(s,t,p)=\cases (J_0(p),0), \; &s\le 0,\\
 (\bar J(s,t,p), \beta(s)\tilde H(t,p)), \;&  0<s\le R, \\
 (\bar J(R+1-s,t,p), \beta(R+1-s)\tilde H(t,p)), \;& R<s\le R+1,\\
 (J_0(p),0), & s>R+1. \endcases$$
 Associated to $(\tilde J_R,\tilde H_R)$ we have the
 Cauchy-Riemann type operator $\bar \partial_R$ for $u:\r\times [0,1]\to
 M$ satisfying the boundary conditions $u(\cdot, 0)\in L$ and $u(\cdot, 1)\in
 L$, and consider the following equation,
 $$\bar \partial_Ru(s,t):=\partial_su+\tilde J_R(s,t,u)(\partial_tu
 +X_{\tilde H_R}(s,t,u))=0. \tag 4.2$$
We note that for $1\le s\le R$, (4.2) describes  the `` negative
gradient flow" for the action functional $\Cal A_H$, i.e., it
satisfies
$$\bar\partial_{J,H}u(s,t):=  {\partial_s  {u}} +\tilde {J}_t\left( {\partial_t   u}
+X_{\tilde{H}}(  {u}(s,t))\right)=0.\tag 4.3$$
 The energy of $u:\r\times [0,1]\to M$ associated to $\tilde J_R$
 is defined by
 $$E_R(u)=\int_{-\infty}^{+\infty}\int_0^1 |\partial_su|_{\tilde J_R}^2\,dsdt.$$
Since a solution $u$ of (4.2) restrict to $(-\infty, 0)\times
[0,1]$ or $(R+1, +\infty)\times [0,1]$ is $J_0$-holomorphic,
finite energy $E_R(u)<\infty$ implies by the boundary removal of
singularities (see [Oh1]) that $u$ can be extended over the disc
carrying the conformal structure from $\r\times [0,1]$,
$$\tilde D=\{-\infty\}\cup  (-\infty,+\infty)\times [0,1]\cup \{+\infty\}.$$
Thus we can identify $\tilde D$ with the standard disc $(D,i)$,
and for every finite energy solution $u$ of (4.2), the homotopy
class $[u]\in \pi(M,L)$ is well defined. We define the adapted
solution spaces associated with $R$ by
$$\Cal M^0(R)=\{u\in C^{\infty}(\r\times [0,1]\to M)\;|\;\bar\partial_R(u)=0,\;E_R(u)<\infty, \;
[u]=0\in \pi_2(M,L)\}.$$
 For an adapted solution $u$, the following result give an
 estimate of the energy of $u$.
  \proclaim{Corollary 4.1} Every solution $u\in \Cal M^0(R)$
 satisfies the energy estimate
 $$0\le E_R(u)\le \|H\|,\;\; \forall R\ge 1. \tag 4.4 $$
 Moreover,  there exists an $l\in \r$ such that
 $$\Cal A_H(u(\varrho,\cdot))\in [l,l+\|H\|],\;\; \forall \varrho\in [1,R]. \tag 4.5$$
 \endproclaim
\demo{Proof} These results are taken from [MS1] (Corollary 4.2)
for the case of  fixed points of Hamiltonian diffeomorphism. The
proof is the same. We give the proof here for the readers'
convenience. For $u\in \Cal M^0(R)$ and $1\le \sigma\le \sigma'\le
R$, there holds
$$0\le E(u^-_{\sigma})\le \Cal A_H(u(\sigma,\cdot), [u^-_{\sigma}])-\int^1_0\inf_{p\in M} H(t,p)\,dt, \tag 4.6$$
$$0\le E(u^+_{\sigma})\le -\Cal A_H(u(\sigma,\cdot), [u^-_{\sigma}])+\int^1_0\sup_{p\in M} H(t,p)\,dt, \tag 4.7$$
$$0\le E(u^-_{\sigma'})-E(u^-_{\sigma})=\Cal A_H(u(\sigma',\cdot), [u^-_{\sigma'}])-
\Cal A_H(u(\sigma,\cdot), [u^-_{\sigma}]). \tag 4.8 $$
 Here $u^-_{\sigma}$ is the restriction of $u$ to $D^-_{\sigma}:=\{-\infty\}\cup(-\infty,\sigma)\times
 [0,1]$ and $u^+_{\sigma}$ is the restriction of $u$ to $D^+_{\sigma}:=(\sigma, +\infty)\times
 [0,1] \cup \{+\infty\}$. (4.6) follows by
 $$\aligned E(u^-_{\sigma})&= \iint_{D^-_{\sigma}}\omega(\partial_su,\tilde J_R\partial_s u)\,dsdt=
  \iint_{D^-_{\sigma}}\omega(\partial_su, \partial_t u+\beta X_{\tilde H})\,dsdt\\
  &=  \iint_{D^-_{\sigma}}u^*\omega+\int^1_0  \tilde H(t,u(\sigma,t))\,dt-\int^{\sigma}_{-\infty}
  \beta'(s)\,ds\int^1_0\tilde H(t,u(s,t))\,dt. \endaligned$$
  Thus there holds
  $$ \Cal A_H(u(\sigma,\cdot),[u^-_{\sigma}])-\int^1_0\sup_{p\in M}H(t,p)\,dt\le E(u^-_{\sigma})\le
  \Cal A_H(u(\sigma,\cdot),[u^-_{\sigma}])-\int^1_0\inf_{p\in M}H(t,p)\,dt. \tag 4.9$$
Using $\Cal A_H(u(\sigma,\cdot),[u^+_{\sigma}])=\omega([u])-\Cal
A_H(u(\sigma,\cdot),[u^-_{\sigma}])$ and $\omega([u])=0$, we get
(4.7) analogously. (4.8) is obvious. (4.4) follows from (4.6) and
(4.7). (4.5) follows from (4.9) and the fact $E(u^-_{\sigma})\le
E_R(u) $.
 \hfill $\blacksquare$ \enddemo

For the modified pair $(\tilde J_R,\tilde H_R)$, as  in section 3,
we choose an auxiliary Morse function $f$ and $1$-parameter
families metrics $g^j_s$ on $L$, $j=0,1,\cdots,k$. For any
$k+1$-tuple $(y_0,\cdots,y_k)\in (\text{Crit}f)^{k+1}$, we define
the moduli space
$$\aligned & \Cal M^0_{y_0;y_1,\cdots,y_k}(J,H,f,(g^j_s))\\
&= \{ (u,\gamma_0,\cdots,\gamma_k)\in \Cal M^0((k+1)R)\times
W^u_{g^0}(y_0)\times W^s_{g^1}(y_1)\times\cdots\times
W^s_{g^k}(y_k) |\\
& \quad\quad\quad \quad\quad\quad\quad\quad\quad\;
u(-\infty)=\gamma_0(0),\;\; u(jR,0)=\gamma_j(0), \;\; j=1,\cdots,k
\}.\endaligned \tag 4.10$$
 Here we remind that we have replace the disc $D$ by the disc
 $\tilde D=\{-\infty\}\cup  (-\infty,+\infty)\times [0,1]\cup
 \{+\infty\}$ with the standard complex structure $i$, and $z_0=-\infty$,
 $z_j=(jR,0)$.

 An immediate consequence of Theorem 3.1 is
 \proclaim {Corollary 4.2} Let $(M, \omega)$ be a closed symplectic
 manifold, $L$ be its closed rational Lagrangian submanifold with the constant $\sigma(L)$ as defined in section
 2. The Hamiltonian $H: [0,1]\times M\to \r$ satisfies
 $\|H\|<\sigma(L)$. Given homogeneous cohomology classes $\alpha_1,\cdots,\alpha_k\in
 H^*(L)$ with nontrivial cup product $\alpha_0=\alpha_1\cup\cdots\cup\alpha_k\in
 H^*(L)$, there exist critical points $y_0,y_1,\cdots,y_k\in
 \text{Crit}f$satisfying
 $$\mu(y_0)=\text{deg}\, \alpha_0,\;\mu(y_j)=\text{deg}\, \alpha_j,\;\;j=1,\cdots,k$$
 such that the solution space $\Cal
 M^0_{y_0;y_1,\cdots,y_k}(J,H,f,(g^j_s))$is nonempty.
 \endproclaim
From this existence result for finite energy solutions of (4.2),
we will deduce the asserted estimate for the number of critical
values for the action functional $\Cal A_H$ by considering $R\to
\infty $.

We now consider the broken flow trajectories. Let us recall the
pair $(\tilde J, \tilde H)$ and the Cauchy-Riemann type equation
from (2.1) with $L$ boundary conditions
$$\aligned &(\bar\partial_{J,H}u)(s,t)=\partial_su+\tilde J(t,u)(\partial_tu+X_{\tilde H}(u))=0,\\
& u(s,0)\in L,\;\; \forall s\in (-\infty,+\infty)\\
         & u(s,1)\in L, \;\;\forall s\in (-\infty,+\infty). \endaligned\tag 4.11$$
 and the Hamiltonian systems with the $L$ boundary
conditions from (2.7)
$$\cases &\dot x(t)=-X_{\tilde H}(x(t))\\& x(0)\in L,\;\;x(1)\in L.\endcases \tag 4.12$$
The set of solutions of (4.12) is $\text {1-1}$ correspondent with
the set of the intersection points $L\cap \varphi^1(L)$. We denote
the  set of solutions of (4.12) by $\Cal S_{L}(H)$.

 \proclaim{Proposition 4.3} If the number of the above solution set  $\sharp \Cal S_{L}(H)<\infty$, then
 there exists a unique limit $x\in \Cal S_{L}(H)$ for every
 solution of (4.11) restrict in the half area with the same
 boundary condition
  $$\aligned &(\bar\partial_{J,H}u)(s,t)=\partial_su+\tilde J(t,u)(\partial_tu+X_{\tilde H}(u))=0,\\
& u(s,0)\in L,\;\; \forall s\in [0,+\infty)\\
         & u(s,1)\in L, \;\;\forall s\in [0,+\infty)\\
         & E(u)<\infty. \endaligned\tag 4.13$$
         that is, $u(s,\cdot)\to x$ uniformly in
         $C^{\infty}([0,1],M)$ as $s\to \infty$.
\endproclaim
\demo{Proof} This proposition is adapted from Proposition 4.4 of
[MS1] and the proof is standard as given in [MS1]. We consider the
reparametrized solution $u_n=u(\cdot+s_n,\cdot)$ for $s_n\to
\infty$, we have $E(u_n|_{[-\sigma,\sigma]})\to 0$ for all
$\sigma>0$ due to the finite energy assumption. Hence for a
suitable subsequence $u_{n_k}$ converges in $C^{\infty}_{loc}$ and
the limit is a translation invariant solution of
$\partial_{J,H}u=0$ with the mentioned boundary conditions over
$\r\times [0,1]$, that is constant in $s$ and therefore an $x\in
\Cal S_{L}(H)$. Given two sequences $s_n, s'_n\to \infty$ with
$u(s_n)\to x$ and $u(s'_n)\to x'$ the finiteness of $\Cal
S_{L}(H)$ implies $x=x'$. Otherwise, one can assume that
$s'_n-s_n\to \infty$ and find, after choosing suitable
subsequence, a sequence $s_n<\tilde s_n<s'_n$ such that without
loss of generality $u(\tilde s_n)\to \tilde x$ with $x\ne \tilde
x$ and $x'\ne \tilde x$. Repeating this argument finitely many
times leads to a contradiction.

 \hfill $\blacksquare$ \enddemo

 Without loss of generality we can assume that $\sharp \Cal
 S_{L}(H)<\infty$. Hence for a solution of (4.11) with finite
 energy, there exist $x,x'\in  \Cal S_{L}(H)$ such that
 $$\lim_{s\to -\infty}u(s)=x,\;\; \lim_{s\to  \infty}u(s)=x'.$$
We define the following connected trajectory spaces for $x,x'\in
\Cal S_{L}(H)$
$$\Cal M_{x,x'}(J,H)=\{u:\r\times [0,1]\to M|\; u\,\, \text{solves}\,\, (4.11),
\lim_{s\to -\infty}u(s)=x,\;\; \lim_{s\to  \infty}u(s)=x'\}. $$
Similarly we define disk type solution spaces for the structure
$\bar J$ and $\beta$ from above
$$\aligned\Cal M_x^{\mp}(\bar J, H)=\{u:\r\times [0,1]\to M|\;\;&\partial_su+
\bar J(\pm s, t,u)(\partial_tu+\beta(\pm s)X_{\tilde H}(t,u))=0\\
&u(s,0)\in L,\;\;  u(s,1)\in L, \;\;\forall s\in \r\\
         & E(u)<\infty,\;\;u(\pm \infty)=x\}.
\endaligned
$$
An element of $\Cal M_x^{\mp}(\bar J, H)$ is a map which is
pseudo-holomorphic in an area containing infinity (the singularity
at infinity can be removed) and is a solution of (4.13) in another
area containing infinity with $x$ as its limit.

We denote the spaces of so-called broken solutions by
$$\aligned\tilde {\Cal M^0}(\bar J,H)=\{&(u_-,u_1,\cdots,u_k,u_+)\\
&\in \Cal M_{x_{0}}^{-}(\bar J,H)\times \Cal
M_{x_0,x_1}(J,H)\times\cdots\times\Cal
M_{x_{k-1},x_k}(J,H)\times\Cal M_{x_{k}}^{+}(\bar J,H)|\\
&x_0,\cdots,x_k\in\Cal S_{L}(H),\;\;k\ge
0,\;[u_-\#u_1\#\cdots\#u_+]=0\in \pi_2(M,L)\},
\endaligned$$
where $\#$ is the obvious gluing operation.

 Considering the solution spaces $\Cal M^0(R_n)$ for $R_n\to
 \infty$, we say that a sequence $u_n\in \Cal M^0(R_n)$ {\it converges
 weakly} to a broken solution
  $$u_n\rightharpoonup (v_0,v_1,\cdots,v_{k},v_{k+1})\in \tilde {\Cal M^0}(\bar J,H) $$
if there are sequences $\{\sigma_{i,n}\}_{n\in \n}\subset \r$,
$i=0,\cdots,k+1$, such that the reparametrized maps
$u_n(\cdot+\sigma_{i,n},\cdot)$ converge uniformly on compact
subsets with all derivatives to $v_i$,
$$u_n(\cdot+\sigma_{i,n},\cdot)\to v_i\;\text{in}\;C^{\infty}_{loc}(\r\times[0,1],M). $$
Clearly, this requires that $\sigma_{0,n}=0$ and
$\sigma_{k+1,n}=R_n+1$ for all $n\in \n$. The following result is
analogous to Gromov's result about the minimal energy of
$J$-holomorphic discs, [G].
 \proclaim{Lemma 4.4} Given a pair $(J,H)$ with   $\sharp \Cal
 S_{L}(H)<\infty$, there exists a lower bound $\hbar(J,H)>0$ for
 the energy of all non-stationary finite energy trajectories,
 that is,
 $$\bar\partial_{J,H}u=0,\;u(s,0)\in L,\;u(s,1)\in L, \;\text{and}\;
 \partial_su\ne 0\;\text{imply}\;E(u)\ge \hbar(J,H). $$

\endproclaim
\demo{Proof} We follow the ideas of [HS] to prove the result. For
the case $H\equiv 0$, $u$ can be extended to a $J$-holomorphic
disc. The result follows from the Gromov compactness. In fact, if
there is a sequence of $J$-holomorphic discs $u_n$ with energy
$E(u_n)\to 0$, then by Gromov compactness, $u_n$ weakly converges
to a cusp curve with positive energy,  a contradiction. If $H\neq
0$, assume that there is a sequence of solution $u_n$ with $0\ne
E(u_n)\to 0$. We prove that
 $\partial_s u$ converges to zero uniformly in $\r\times [0,1]$ as
 $n$ tends to $\infty$. Otherwise there would exist a sequence
 $(s_n,t_n)$such that $|\partial_s u(s_n,t_n)|\ge \delta >0$.
If $s_n$ is bounded, we can assume $s_n\to 0$ without loss of
generality. Since $E(u_n)$
 converges to zero no bubbling can occur and hence a subsequence
 of $u_n$ converges with its derivatives uniformly on compact sets
 to a solution $u:\r\times [0,1]\to M$ with mentioned boundary
 conditions, $\partial_su(0,t^*)\ge \delta$ and $E(u)=0$. But the
 latter implies that $u(s,t)\equiv x(t)$ in contradiction to the
 former. If $s_n$ is non-bounded, then we can assume $s_n\to
 \infty$. We consider $v_n(s,t)=u(s+s_n, t)$ as in the proof of
 Proposition 4.3, then by the finiteness condition: $\sharp \Cal
 S_{L}(H)<\infty$, we can get $v_n\to v$ with $|\partial_s v(0,t^*)|\ge \delta$
 and $E(v)=0$, it is still a contradiction.
 \hfill $\blacksquare$ \enddemo

 We denote the broken trajectory space by
$$\aligned \Bar {\Cal M}_{x,y}(J,H)=\{&\bold u=(u_1,\cdots,u_r)\,|\;u_i\in\Cal M_{x_{i-1},x_i}\\
& i=1,\cdots,r,\;x_0=x,\;x_r=y,\; r\in\n\}.\endaligned $$
 It is the space of broken trajectories started from $x\in \Cal S_L(H)$
 and ended at $y\in \Cal S_L(H)$. The energy of a broken
 trajectory $\bold u=(u_1,\cdots, u_r)$ satisfies
 $$E(\bold u)=\sum_{i=1}^{r}E(u_i). $$
 If $\bold u\in \bar{\Cal M}_{x,x}$, then $[\bold u]\in\pi_2(M,L)$
 is well defined and $\omega([\bold u])=E(\bold u)\ne 0$, the latter follows from the fact that the
 start point is just the end point, so there holds
 $$ \sum_{i=1}^r\int_{-\infty}^{\infty}
 \omega(\partial_t u_i, JX_{\tilde H}(u))\,dt=0. $$
 Thus if $\bold u\ne x$, then
 $E(\bold u)\ne 0$, it implies that $\omega([\bold u])=E(\bold u)\ge
 \sigma(L)$.

 We define
  $$\sigma_0(\omega,H,J)=\inf\{E(\bold u)\,|\,\bold u\in \bar \Cal M_{x,x}(J,H),\; \bold u\neq x,
  \; x\in \Cal  S_{L}(H)\}.$$

\proclaim{Theorem 4.5} Let $\sharp \Cal  S_{L}(H)<\infty$ and
$u_n\in \Cal M^0(R_n)$ be a sequence of solution with $R_n\to
\infty$ and uniformly bounded gradient $\nabla u_n$. Then there
exists a subsequence $\{\u_{n_k}\}$ converging weakly to a broken
solution
$$u_{n_k}\rightharpoonup (v_-,v_1,\cdots,v_N,v_+)\in \Cal M^0(\bar J, H). $$
\endproclaim
  \demo{Proof} This result is similar to Theorem 4.5 of [MS1].
  Elliptic bootstrapping implies $C^{\infty}_{loc}$-convergence
  for subsequences of $\{u_n(\cdot+s_n)\}$ for any shifting
  sequences $\{s_n\}$, $s_n\to \infty$. Assume that we have
  already shifting sequences $\{s_n\}$ and $\{\bar s_n\}$ such
  that $s_n-\bar s_n\to \infty$ and $u_n(\cdot+s_n)\to v,\; u_n(\cdot+\bar s_n)\to w$
  in $C^{\infty}_{loc}$ with $v\in \Cal M_{x,y}(J,H)$ and $w\in \Cal
  M_{y',z}(J,H)$, we use the analogous argument as in the proof of
  Proposition 4.3. We show that either $y=y'$ or that modulo
  choosing a subsequence we find a sequence $\tau_n\to \infty$
  such that $s_n<\tau_n<\bar s_n$ and $u_n(\cdot+\tau_n)\to \bar w\in \Cal
  M_{y,y'}(J,H)$. This requires lifting to the covering
  $\Omega_0(L)$ where the function $\Cal A_H$ is real-valued and
  the energy of $u\in \Cal M_{x,y}(J,H)$ is given by $E(u)=\Cal A_H(\bold y)-\Cal A_H(\bold
  x)$,   where $\bold x$ is the lifting of $x$ in $\Omega_0(L)$.
  From the total energy bound by $\|H\|$ from Corollary 4.1 and
  the minimal energy $\hbar(J,H)>0$ for non-stationary
  trajectories from Lemma 4.4, it follows that only finite number
  of $\bar \bold y\in \tilde \Cal  S_{L}(H)$, the lifting of $\Cal
  S_{L}(H)$ in  $\Omega_0(L)$, can occur between $\bold y$ and $\bold y'$. It
  remains to show that $\Cal A_H(\bold y)=\Cal A_H(\bold y')$
  implies $  y=  y'$. This follows from the following
  result.

\proclaim{Lemma 4.6} Let $\sharp \Cal  S_{L}(H)<\infty$, there
exists a $\gamma>0$ such that for every neighbourhood $W$ of $\Cal
S_{L}(H)$ in $C^{\infty}([0,1],M)$ there exists a number $h=h(M)$
with the following properties:

If $u: (r,R)\times [0,1]\to M$ for $-\infty\le r<R\le \infty$
solves
$$\aligned &\bar\partial_{J,H}u=0,\;\;u(\cdot,0)\in L, \;\; u(\cdot,1)\in L,\;\;[u(\frac{r+R}{2},\cdot)]=0\in \pi_1(M,L),\\
&E(u)\le \gamma\;\;\text{and}\;\;R-r>2h,
\endaligned  \tag 4.14$$
then $u(s)\in W$ for all $s\in (r+h,R-h)$. Moreover, given $k_0\in
\n,\;\epsilon>0 $, there exists $h=h(k_0,\epsilon)$ such that
solutions of (4.14) view as a mappings into $M\subset \r^N$
satisfy
$$|D^{\alpha}(u(s,t)-x(t))|\le \epsilon, \forall\, (s.t)\in (r+h,R-h)\times [0,1],\;|\alpha|\le k_0 $$
for a suitable $x\in \Cal  S_{L}(H)$.
\endproclaim
  \demo{Proof} We prove indirectly the second assertion. Assume
  that given any $\gamma>0$ there exist $k(\gamma)\in \n,\;\epsilon(\gamma)>0,\; h_n\to \infty,\; r_n<R_n$
  with $R_n-r_n\ge 2h_n$ and $u_n:(r_n,R_n)\times [0,1]\to M$
  satisfying the boundary condition as in (4.14) and
  $\bar\partial_{J,H}u_n=0$, $\int^{R_n}_{r_n}\int^1_0\,|\partial_su_n|^2\,dsdt\le \gamma$
  such that there exist $(s_n,t_n)\in (r_n,R_n)\times [0,1]$ and
  $\alpha\le k$ with
  $$|D^{\alpha}(u_n(s_n,t_n)-x(t_n))|>\epsilon $$
for all $n\in\n$ and $x\in \Cal  S_{L}(H)$. Reparametrizing $u_n$
so that $v_n(s,t)=u_n(s+s_n,t)$ solves $\bar
\partial_{J,H}v_n=0$ with
$$\int^{h_n}_{-h_n}\int^1_0 |\partial_sv_n|^2\,dsdt<\gamma\;\;\text{and}\;\;
|D^{\alpha}(v_n(0,t_n)-x(t_n))|>\epsilon.$$
 Without loss of generality we can replace $t_n$ by some $t_0$.
Choosing $\gamma>0$ small enough by Gromov's theorem about the
minimal energy of pseudoholomorphic spheres or holomorphic discs
with $L$-boundary condition (see the proof of Lemma 4.4), there
exists a number $c>0$ such that
 $$|\nabla v_n(s,t)|\le c\;\;\forall\, (s,t)\in [-\frac 34h_n,\frac34h_n]\times [0,1],\;\;n\in\n. $$
 Otherwise, we would obtain a pseudoholomorphic sphere or disc
 bubbling off with energy less than $\gamma$. Thus , choosing a
 suitable subsequence, without loss of generality denoted again by
 $n\in \n$, we obtain uniform convergence on compact subsets, $v_n\to (v:\r\times [0,1]\to
 M)$ in $C^{\infty}_{loc}$
 with
 $$\aligned & \bar\partial_{J,H}v=0,\;\;\int^{\infty}_{-\infty}\int^1_0|\partial_sv|^2\le \gamma\;\;\text{and}\\
 &|D^{\alpha}(v(0,\cdot)-x(\cdot))|_{L^{\infty}([0,1])} >\epsilon,\;\; v(\cdot,0)
 \in L,\;\;v(\cdot,1)\in L.\endaligned$$
 But Lemma 4.4 implies for $\gamma<\hbar(J,H)$ that
 $\partial_sv=0$, i.e. $v(0)\in \Cal  S_{L}(H)$ providng the
 contradiction.
  \hfill $\blacksquare$ \enddemo
  This also concludes the proof of Theorem 4.5 because $\Cal A_H(\bold y)=\Cal A_H(\bold y')$
  implies that we can find sequences $s_n$ and $s_n'$ such that $u_n(s_n)\to
  y$, $u_n(s_n')\to y'$ and $0<s_n-s_n'$  with $E(u_n|_{[s_n,s_n']})\to
  0$. Consequently, Lemma 4.6 yields $y=y'$. \hfill $\blacksquare$ \enddemo
{\bf Remark 4.7.} In our case we have $\|H\|<\sigma(L)$, and
$E(u_n)\le \|H\|$ by Corollary 4.1 for $u_n\in \Cal M^0(R_n)$,
bubbling-off cannot occur, thus the gradient of $u_n$ is uniformly
bounded.

 We
denote the covering space of $\Cal S_L(H)$ in the sense of section
2 by $\tilde \Cal S_L(H)$, i.e., any element $\bold x=(x,u_x)\in
\tilde \Cal S_L(H)$ is a critical point of $\Cal A_H$ in the space
$\Omega_0(L)$ and $x$ is a solution of $\dot x=-X_{\tilde H}(t,x)$
with the boundary conditions $x(0)\in L$ and $x(1)\in L$. It
implies $x(0)\in L\cap \varphi^1(L)$, see (2.7). The space carries
a partial ordering with respect to the gradient flow of $\Cal
A_H$.
 \proclaim{Definition 4.8} Given a pair $\bold x,\;\bold x'\in \tilde \Cal
 S_L(H)$, we say $\bold x\le \bold x'$ if there exist connecting
 broken flow trajectories $\bar{\Cal M}_{x,x'}(J,H)\ne \varnothing$.
 Given a Morse-Smale pair $(f,g)$, we say that $\bold x\ll \bold
 x'$if there exist $u\in \Cal M_{x,x'}(J,H)$ and $y\in
 \text{Crit}f$ with $\mu(y)\ge 1$ such that $u(0,0)\in W^s_g(y)$.
 \endproclaim

If $\sharp\{L\cap \varphi^1(L)\}<\infty$, then for a generic
choice of Morse function $f: L\to \r$  and  Riemannian metric $g$
on $L$ , there holds
$$\bigcup_{\mu(y)\ge 1}W^s_g(y)\cap L\cap \varphi^1(L)=\varnothing. \tag 4.15$$
This can be prove by standard transversal analysis (see [MS1]).
Thus if choose $(f,g)$ satisfying (4.15), then for $\bold x\ll
\bold x'$ we have $\bold x\ne\bold x'$ thus $\bold x<\bold x'$ and
in particular $\Cal A_H(\bold x)<\Cal A_H(\bold x')$. The latter
can be seen from  the proof of Theorem 4.5. By this observation we
have the following result.
 \proclaim{Corollary 4.9} Let $k\in \n$ and $(f,g^i_s)$, $i=1,\cdots,k$ satisfy
 condition  (4.15) with respect to $H$ satisfying $\Cal
 S_L(H)<\infty$. Given a sequence
  $$u_n\in \Cal M^0_{y_0;y_1,\cdots,y_k}((k+1)R_n),\;R_n\to \infty $$
  with $y_i\in \text{Crit}f$, $\mu(y_i)\ge 1$ for
  $i=0,1,\cdots,k$, weakly converging to a broken trajectory,
  there exist solutions $\bold x_1,\cdots,\bold x_N\in \tilde \Cal
  S_L(H)$ satisfying $\bold x_1\le\cdots\le\bold x_N$ and
  $1\le n_1<m_1\le n_2<m_2\le\cdots\le n_k<m_k\le N$ such that $\bold x_{n_i}\ll \bold x_{m_i}$
  for $i=1,\cdots,k$. In particular, there exists an $l\in \r$ such that
  $$l\le \Cal A_H(\bold x_{n_1})<\cdots<\Cal A_H(\bold x_{n_k})<
  \Cal A_H(\bold x_{m_k})\le l+\|H\|. $$
  \endproclaim
\demo{Proof} By assumption, the sequence $u_n\in\Cal
M^0((k+1)R_n)$ satisfies
 $$u_n(jR_n,0)\in W^s_g(y_j),\;\;j=1,\cdots,k. $$
 Moreover, if $u_n$ converges weakly to a broken solution
 $$(v_-,v_1,\cdots,v_N,v_+)\in \tilde \Cal M^0(\bar J,H) $$
 we have reparametrization sequences $\{\sigma_{i,n}\}_{n\in\n}$
 for $i=1,\cdots,N$ such that $u_n(\cdot+\sigma_{i,n},\cdot)\to
 v_i$ in $C^{\infty}_{loc}$ and $u_n\to v_-$, $u_n(\cdot-(k+1)R_n-1,\cdot)\to
 v_+$. Considering the shifted solutions $u_{n,j}=u_n(\cdot-jR_n,\cdot)$,
 we thus obtain after choosing a suitable subsequence
 $C^{\infty}_{loc}$-convergence $u_{n,j}\to w_j\in \Cal M_{ x_j,  x_j'}(J,H)$ for some
 $\bold x_j,\bold x_j'\in \tilde\Cal S_L(H)$, $j=1,\cdots,k$. By
 definition, we have $\bold x_j\ll\bold x_j'$ and the assumption
 of weak convergence implies the order
 $$\bold x_1\ll \bold x'_1\le \bold x_2\ll\bold x'_2\le \cdots\le \bold x_k\ll\bold x'_k. $$
We now can prove the main result of this paper.\hfill
$\blacksquare$ \enddemo
 \proclaim{Theorem 4.10} Let $(M,\omega)$ be a closed symplectic manifold, and
 $L$ be its closed Lagrangian submanifold satisfying the rational
 condition $\omega(\pi_2(M,L))=\sigma(L)\cdot \z,\;\sigma(L)>0$.
 $\varphi=\varphi^1$ is a Hamiltonian automorphism of $(M,\omega)$
 generated by the Hamiltonian $H:[0,1]\times M\to \r$ with
 $\|H\|<\sigma(L)$. Then the cup-length estimate of the Lagrangian intersection holds
 $$\sharp\{L\cap \varphi(L)\}\ge \text{cl}(L).$$
 \endproclaim
\demo{Proof} By the assumption $\|H\|<\sigma(L)$, for a generic
almost complex structure $J$ compatible with the symplectic
structure $\omega$, let $k+1=\text{cl}(L)$, then by Corollary 4.2
we find solutions $u_n\in \Cal M^0_{y_0;y_1,\cdots,y_k}((k+1)R_n)$
for some sequence $R_n\to \infty$ and $y_i\in \text{Crit}f$ where
$(f,g^i)$ satisfy (4.15). By Corollary 4.9, Theorem 4.5 and Remark
4.7, there are $k+1$ critical points $\bold x_i\in \tilde \Cal
S_L(H)$ for $\Cal A_H$ on $\Omega_0(L)$ defined in section 2 such
that
$$l\le \Cal A_H(\bold x_1)<\cdots<\Cal A_H(\bold x_{k+1})\le l+\|H\| $$
for some $l\in \r$. Due to the assumption $\|H\|<\sigma(L)$ again,
there is no broken trajectory of flow started from some solution
$x\in \Cal S_L(H)$ and ended at the same solution. In fact, the
energy of this mentioned broken trajectory should be not less than
the number $\sigma(L)$, but on the other hand side, this energy
should be not more than $\|H\|$ since $E(u_n)\le \|H\|$. Namely,
the $k+1$ critical points $\bold x_i$ project to $k+1$ different
solutions $x_i\in \Cal S_L(H)$.
 \hfill $\blacksquare$ \enddemo

 {\bf Remark 4.12.} {\bf Remark.} As  in [Ch1-Ch2], the
symplectic manifold can be
 more generally a tame symplectic manifold, since the tameness condition allows us to deal with $M$
as if it is compact, all the techniques are the same as in the
compact case if we only consider the compactly supported
Hamiltonian $H$. We recall that $(M,\omega)$ is tame if  there
exists an almost complex structure $J$ on $M$ such that
$g(\cdot,\cdot)=\omega(\cdot,J\cdot)$ is a Riemannian metric on
$M$ satisfying the following conditions:

 (i) Riemannian manifold $(M,g)$ is complete,

 (ii) the sectional curvature of $g$ is bounded,

 (iii) the injectivity radius of $g$ is bounded away from zero.

 Let $J$ be an almost complex structure on $M$ such that
 $(M,\omega,J)$ is a tame almost K\"{a}hler manifold, denote by $\Cal
 J$ the  space of such stuctures. Let $\sigma_S(M,J)$ denote the
 minimal area of a $J$-holomorphic sphere in $M$, and
 $\sigma_D(M,L,J)$ denote the minimal area of a $J$-holomorphic
 disc in $M$ with boundary on $L$. These numbers may equal
 infinity if there are no such $J$-holomorphic curves. Otherwise,
 minimals are achieves due to the Gromov compactness theorem (see
 [G])and are clearly positive. Let
 $$\sigma(M,L,J)=\min(\sigma_S(M,J),\sigma_D(M,L,J)) $$
We remind the number $\sigma_0(\omega,H,J)$ is defined just before
Theorem 4.5. The following result does not require that $L$ is
rational Lagrangian submanifold of $M$.
 \proclaim{Theorem 4.13} If $\|H\|<\min(\sigma_0(\omega,H,J),
 \sigma(M,L,J))$, then the standard cup-length estimate is valid
 $$\sharp (L\cap \varphi(L))\ge \text{\rm cl}(M). $$
 \endproclaim

\demo{Proof} The proof is the same as in the proof of Theorem
4.11. With the condition $\|H\|<\min(\sigma_0(\omega,H,J),
 \sigma(M,L,J))$, the bubbling-off can not occur, we can also guarantee that the
 different critical points $\bold x_i\in \tilde \Cal S_L(H)$ can be
 project to different $x_i\in \Cal S_L(H)$ as done in the proof of
 Theorem 4.11.
 \hfill $\blacksquare$ \enddemo

\enddocument